\documentclass[proceedings,submission]{dmtcs}
\usepackage{amssymb, amscd, amsmath, amsfonts}
\usepackage[dvipsnames]{xcolor}
\usepackage{verbatim}

\newcommand{\idiot}[1]{\vspace{5 mm}\par \noindent
\marginpar{\textsc{For longer version}}
\framebox{\begin{minipage}[c]{.99 \textwidth}
#1 \end{minipage}}\vspace{5 mm}\par}
 
\newcommand{\todone}[1]{\vspace{5 mm}\par \noindent
\marginpar{\textsc{DONE!}}
\framebox{\begin{minipage}[c]{.99 \textwidth}
\tt #1 \end{minipage}}\vspace{5 mm}\par}
 
\renewcommand{\todone}[1]{}
 
\renewcommand{\idiot}[1]{}
 
\newdimen\squaresize \squaresize=12pt
\newdimen\thickness \thickness=0.4pt
 
\def\square#1{\hbox{\vrule width \thickness
   \vbox to \squaresize{\hrule height \thickness\vss
      \hbox to \squaresize{\hss#1\hss}
    \vss\hrule height\thickness}
\unskip\vrule width \thickness}
\kern-\thickness}
 
\def\vsquare#1{\vbox{\square{$#1$}}\kern-\thickness}

\def\thisbox#1{\kern-.09ex\fbox{#1}}
\def\downbox#1{\lower1.200em\hbox{#1}}

\def\boldentry#1{\textcolor{Red}{\textbf{#1}}}

\usepackage{tikz}
\usetikzlibrary{calc}
\usepackage{ifthen}
\newcommand{\tikztableau}[2][scale=0.6,every node/.style={font=\small}]{
    \def\newtableau{#2}
    \begin{array}{c}
    \begin{tikzpicture}[#1]
    \coordinate (x) at (-0.5,0.5);
    \coordinate (y) at (-0.5,0.5);
    \foreach \row in \newtableau {
        \coordinate (x) at ($(x)-(0,1)$);
        \coordinate (y) at (x);
        \foreach \entry in \row {
            \ifthenelse{\equal{\entry}{X}}
               {
                \node (y) at ($(y) + (1,0)$) {};
                \fill[color=gray!10] ($(y)-(0.5,0.5)$) rectangle +(1,1);
                \draw[color=gray] ($(y)-(0.5,0.5)$) rectangle +(1,1);
               }
               {
                \ifthenelse{\equal{\entry}{\boldentry X}}
                   {
                    \node (y) at ($(y) + (1,0)$) {};
                    \fill[color=gray] ($(y)-(0.5,0.5)$) rectangle +(1,1);
                    \draw ($(y)-(0.5,0.5)$) rectangle +(1,1);
                   }
                   {
                    \node (y) at ($(y) + (1,0)$) {\entry};
                    \draw ($(y)-(0.5,0.5)$) rectangle +(1,1);
                   }
               }
            }
        }
    \end{tikzpicture}
    \end{array}}

\newcommand{\tikztableausmall}[1]{\tikztableau[scale=0.42,every node/.style={font=\small}]{#1}}

\usepackage[latin1]{inputenc}
\usepackage{subfigure}

\def\ZZ{\mathbb{Z}}

\def\BB{\mathbb{B}}

\def\B{{\bf B}}

\def\sym{\operatorname{\mathsf{Sym}}}

\def\Qsym{\operatorname{\mathsf{QSym}}}

\def \fS{{\mathfrak S}}
\def \HH{{H}}
\def\Nsym{\operatorname{\mathsf{NSym}}}

\address{\addressmark{1} Fields Institute, Toronto, ON, Canada\\
\addressmark{2}Laboratoire de combinatoire et d'informatique math\'ematique, Universit\'e du Qu\'ebec \`a Montr\'eal, Montr\'eal, QC, Canada\\
\addressmark{3}
York University, Toronto, ON, Canada
}

\begin{document}

\newtheorem{Theorem}{Theorem}[section]
\newtheorem{Remark}[Theorem]{Remark}
\newtheorem{Definition}[Theorem]{Definition}
\newtheorem{Lemma}[Theorem]{Lemma}
\newtheorem{Corollary}[Theorem]{Corollary}
\newtheorem{Example}[Theorem]{Example}
\newtheorem{Proposition}[Theorem]{Proposition}
\newtheorem{Conjecture}[Theorem]{Conjecture}

\title[Immaculate basis of the non-commutative symmetric functions]{The immaculate basis of the non-commutative symmetric functions (Extended Abstract)}
\author[C. Berg \and N. Bergeron \and F. Saliola \and L. Serrano \and M. Zabrocki]{Chris Berg\addressmark{2} \and Nantel Bergeron\addressmark{1,3} \and Franco Saliola\addressmark{2} \and Luis Serrano\addressmark{2} \and Mike Zabrocki\addressmark{1,3}}
\date{\today}

\maketitle
\begin{abstract} 
\paragraph{Abstract} We introduce a new basis of the non-commutative symmetric functions whose elements have Schur functions as their commutative images. Dually, we build a basis of the quasi-symmetric functions which expand positively in the fundamental quasi-symmetric functions and decompose Schur functions according to a signed combinatorial formula.

\paragraph{R\'esum\'e.} Nous introduisons une nouvelle base des fonctions sym\'etriques non commutatives dont les images commutatives  sont des fonctions de Schur. Nous construisons la base duale des fonctions quasi-sym\'etriques qui s'expriment de fa\c{c}on positive en fonction de la base  fondamental  et  d\'ecomposer les fonctions de Schur.

\end{abstract}

\setcounter{tocdepth}{3}

\section{Introduction}

The Schur functions $s_\lambda$ are indexed by integer partitions and form an additive basis for the algebra of symmetric functions $\sym$. Schur functions play an important role throughout mathematics, in particular in algebraic geometry (as representatives of Schubert classes for the Grassmannian) and representation theory (they are the characters of the irreducible representations of the general linear group). Another important basis for $\sym$ is the (complete) homogeneous symmetric functions $h_\lambda$.

The algebras of non-commutative symmetric functions $\Nsym$ and quasi-symmetric functions $\Qsym$ are dual Hopf algebras. These algebras have been of great importance to algebraic combinatorics. As seen in \cite{ABS}, they are universal in the category of combinatorial Hopf algebras. They also represent  the Grothendieck rings for the projective and finite dimensional representation theory of the $0$-Hecke algebra \cite{KT}. An important basis for $\Nsym$ is formed by the (complete) homogeneous non-commutative symmetric functions $H_\alpha$, indexed by compositions. The \textit{forgetful map} $\chi: \Nsym \longrightarrow \sym$ maps the homogeneous non-commutative symmetric functions to their symmetric counterparts (see (\ref{def:chi})).

The main goal of this abstract is to define and outline the properties of a new basis, the \textit{immaculate basis} $\fS_\alpha$ of $\Nsym$, which emulates the role of the Schur functions. This new basis projects onto the Schur basis under the forgetful map and it shares many of the same properties and constructions of the classical basis of Schur functions of $\sym$. More specifically:

\textbf{Bernstein operators.} One way to construct Schur functions is by iterating the Bernstein row adding operator, which acts on Schur functions by adding a row to the corresponding Ferrers shape. These operators can be described in an algebraic way, which we deform in order to obtain a non-commutative Bernstein operator. This deformed operator now acts on immaculate functions by adding a row to the corresponding composition. Thus, a repeated iteration of these operators will build the immaculate functions, as in Definition \ref{def:immaculate}.

\textbf{Pieri rule.} The product of a Schur function and a homogeneous symmetric function corresponding to a partition with only one part can be expressed, via the classical Pieri rule, as a multiplicity-free sum over a specific set of Schur functions. More specifically, this sum is over all ways to add a \textit{horizontal strip} to the original shape. In Theorem \ref{thm:Pieri} we show that in a similar way, the product of an immaculate function and a homogeneous non-commutative symmetric function corresponding to a composition with only one part can be expressed as a multiplicity-free sum of immaculate functions. This sum is over all ways to add an analog of a horizontal strip for composition shapes.

\textbf{Immaculate tableaux and the immaculate Kostka matrix.} By iterating the Pieri rule, one can obtain an expansion of the homogeneous symmetric functions in terms of Schur functions, where each coefficient is a \textit{Kostka number}, or number of semistandard Young tableaux of a specified shape and content. In a similar fashion, we introduce \textit{immaculate tableaux}, and by iterating the immaculate Pieri rule, one obtains an expansion of the homogeneous non-commutative functions in terms of the immaculate functions, where each coefficient is the number of immaculate semistandard tableaux of a specified shape and content (Theorem \ref{thm:Hexpansion}).

\textbf{Positive expansion for ribbons.} Another important basis of $\Nsym$ is formed by \textit{ribbon noncommutative functions} $R_\alpha$. In Theorem \ref{thm:Rpositive} we expand the Ribbon functions positively in terms of immaculate functions, indexed by certain descent sets on standard immaculate tableaux.

Moreover, the immaculate basis gives rise to a dual basis in the quasi-symmetric function algebra. The \textit{dual immaculate basis} also shares interesting properties with the Schur basis. In particular, by duality arguments, one is able to express the dual immaculate basis in terms of other known bases of $\Qsym$.

\textbf{Jacobi-Trudi determinant formula.} The Schur functions can be expanded in terms of the homogeneous symmetric functions by the use of the Jacobi-Trudi determinant. By considering a non-commutative version of this determinant, we expand the immaculate functions in terms of the homogeneous non-commutative symmetric functions, thus obtaining a lifting of the Jacobi-Trudi formula in $\Nsym$, as in Theorem \ref{thm:JT}.

\textbf{Generating series of immaculate tableaux and monomial expansion.} The most well known construction for a Schur function is by its expression as a generating series over the set of semistandard Young tableaux, and thus, as a positive sum of monomial (quasi-)symmetric functions. In Theorem \ref{thm:MonomialPositive} we express the dual immaculate functions as a generating series over the set of semistandard immaculate tableaux, and thus, as a positive sum of monomial quasi-symmetric functions.

\textbf{Positive fundamental expansion.} The Schur functions can also be expressed as a positive sum of fundamental quasi-symmetric functions, by considering descents on standard Young tableaux. By a duality argument, in Theorem \ref{thm:FundamentalPositive}, we express the dual immaculate functions as a positive sum of fundamental quasi-symmetric functions, by considering descents on standard immaculate tableaux.

\textbf{Expansion of Schur functions.} In Theorem \ref{thm:decompose}, we show that the Schur functions expand in the dual immaculate basis via signed combinatorics developed in \cite{ELW}.

\textbf{Littlewood-Richardson rule.} In the classical case, the product of two Schur functions can be expressed as a sum of Schur functions, where each coefficient is a \textit{Littlewood-Richardson number}, namely, the number of \textit{Yamanouchi tableaux} of a certain skew shape. Although the product of any two immaculate functions is not in general immaculate positive, we give a combinatorial formula for the coefficients in the product of any immaculate function with an immaculate function corresponding to a partition as the positive sum of immaculate functions, where each coefficient counts the number of \textit{immaculate Yamanouchi tableaux} of a certain skew shape, thus obtaining an analogue of the Littlewood-Richardson rule (Theorem \ref{thm:LR}).

\textbf{Murnaghan-Nakayama rule.} The product of a Schur function and a power sum can be expressed as a sum over Schur functions, over the set of shapes that are obtained by adding a ribbon to the original Ferrers shape. In $\Nsym$, an analogue of the power sums basis $\Psi_\alpha$, was defined in \cite{GKLLRT}. In Theorem \ref{thm:MNrule} we express the product of an immaculate function and a noncommutative power sum $\Psi_n$. 

\textbf{Indecomposable modules.} There exists a collection of indecomposable modules for the $0$-Hecke algebra with the property that the module indexed by the composition $\alpha$ has the dual immaculate function indexed by $\alpha$ as its characteristic. In the interest of space, we will not pursue this below, but refer the reader to \cite{BBSSZ2}.

This text is an extended abstract of the preprints \cite{BBSSZ1}, \cite{BBSSZ2} and \cite{BBSSZ3}, where complete proofs can be found.
\begin{Remark}
Although our basis of $\Nsym$ is similar to the dual basis of quasi-symmetric Schur functions of \cite{HLMvW} (whose properties were developed in \cite{BLvW}), they are in fact different bases.
\end{Remark}
\subsection{Acknowledgments}
This work is supported in part by CRC and NSERC.
It is partially the result of a working session at the Algebraic
Combinatorics Seminar at the Fields Institute with the active
participation of C.~Benedetti, Z.~Chen, H.~Heglin, and D.~Mazur. In addition, the authors would like to thank F. Hivert, J. Huang, J. Remmel,
N. Thi\'ery and M. Yip for lively discussions.

This research was facilitated by computer exploration using the open-source
mathematical software \texttt{Sage}~\cite{sage} and its algebraic
combinatorics features developed by the \texttt{Sage-Combinat}
community \cite{sage-co}.
\section{Background}
\subsection{Compositions and combinatorics} \label{sec:compositions}

A \textit{partition} of a non-negative integer $n$ is a sequence
$\lambda = (\lambda_1, \lambda_2, \dots, \lambda_m)$ of non-negative integers satisfying 
$\lambda_1 \geq \lambda_2 \geq \cdots \geq \lambda_m$, and is denoted 
$\lambda \vdash n$. Partitions are of particular importance to 
algebraic combinatorics; among other things, partitions of $n$ index a basis 
for the symmetric functions of degree $n$, $\sym_n$, and the character ring 
for the representations of the symmetric group. These concepts are intimately 
connected; we assume the reader is well versed in this area (see for instance \cite{Sagan} for background details).

A \textit{composition} of a non-negative integer $n$ is a list 
$\alpha = [\alpha_1, \alpha_2, \dots, \alpha_m]$ of positive 
integers which sum to $n$, written $\alpha \models n$.
The entries $\alpha_i$ of the composition are referred to as the parts
of the composition.  The size of the composition is the sum of the parts
and will be denoted $|\alpha|:=n$.  The length of the composition is the
number of parts and will be denoted $\ell(\alpha):=m$.
In this paper we study 
dual graded Hopf algebras whose bases at level $n$ are indexed by compositions of $n$.

Compositions of $n$ correspond to subsets of $\{1, 2, \dots, n-1\}$. We will follow the convention of identifying $\alpha = [\alpha_1, \alpha_2, \dots, \alpha_m]$ with the subset $D(\alpha) = 
\{\alpha_1, \alpha_1+\alpha_2, \alpha_1+\alpha_2 + \alpha_3, \dots, \alpha_1+\alpha_2+\dots + \alpha_{m-1} \}$. 

If $\alpha$ and
$\beta$ are both compositions of $n$,  say that $\alpha \leq  \beta$ in refinement order if $D(\beta) \subseteq D(\alpha)$. For instance, $[1,1,2,1,3,2,1,4,2] \leq [4,4,2,7]$, since $D([1,1,2,1,3,2,1,4,2]) = \{1,2,4,5,8,10,11,15\}$ and $D([4,4,2,7]) = \{4,8,10\}$.

We  introduce  a  new  notion  which  will  arise  in  our  
Pieri  rule  (Theorem  \ref{thm:Pieri});  we  say  that  $\alpha  \subset_{i}  \beta$  if:
\begin{enumerate}
\item $|\beta| = |\alpha| + i$,
\item $\alpha_j \leq \beta_j$ for all $1 \leq j \leq \ell(\alpha)$,
\item $\ell(\beta) \leq \ell(\alpha) + 1.$
\end{enumerate}
For a composition $\alpha = [\alpha_1, \alpha_2, \dots, \alpha_\ell]$ and a positive integer $m$, we let $[m, \alpha]$ denote the composition $[m, \alpha_1, \alpha_2, \dots, \alpha_\ell]$.
\subsection{Schur functions and creation operators}
We let $h_i$ and $e_i$ denote the complete homogeneous and elementary symmetric functions  of degree $i$ respectively.
We next define a Schur function indexed by an arbitrary sequence of integers.
The resulting family of symmetric functions indexed by partitions $\lambda$ are the usual
Schur basis of the symmetric functions.

\begin{Definition} \label{def:JTformula}
For an arbitrary integer tuple 
$\alpha = (\alpha_1, \alpha_2, \ldots, \alpha_\ell) \in \ZZ^\ell$, we define 
\[ s_\alpha := \det \begin{bmatrix} 
h_{\alpha_1}&h_{\alpha_1+1}&\cdots&h_{\alpha_1+\ell-1}\\
h_{\alpha_2-1}&h_{\alpha_2}&\cdots&h_{\alpha_2+\ell-2}\\
\vdots&\vdots&\ddots&\vdots\\
h_{\alpha_\ell-\ell+1}& h_{\alpha_\ell-\ell+2}&\cdots&h_{\alpha_\ell}\\
\end{bmatrix} = 
\det \left|h_{\alpha_i + j - i}\right|_{1 \leq i,j \leq \ell}\]
where we use the convention that $h_0 = 1$ and $h_{-m} = 0$ for $m>0$.
\end{Definition}

With this definition, switching two adjacent rows of the defining matrix
has the effect of changing the sign of the determinant.  It is also equal to the
Schur function indexed by a different integer tuple:
$$s_{\alpha_1, \alpha_2, \ldots, \alpha_r, \alpha_{r+1}, \ldots, \alpha_\ell} =
-s_{\alpha_1, \alpha_2, \ldots, \alpha_{r+1}-1,\alpha_r+1, \ldots, \alpha_\ell} ~.
$$
\begin{Proposition} \label{prop:schurcomposition}
If $\alpha$ is a composition of $n$ with length equal to $k$, then
$s_\alpha = 0$ if and only if there exists $i,j \in \{ 1,2, \ldots,k\}$ with $i \neq j$ 
such that $\alpha_i -i = \alpha_j-j$.  If $s_\alpha \neq 0$, then
there is a unique permutation $\sigma$ such that $(\alpha_{\sigma_1} +1-\sigma_1, \alpha_{\sigma_2}+2-\sigma_2, \ldots,
\alpha_{\sigma_k}+k-\sigma_k)$ is a partition.  In this case,
$$s_\alpha = (-1)^\sigma s_{\alpha_{\sigma_1} +1-\sigma_1, \alpha_{\sigma_2}+2-\sigma_2, \ldots, \alpha_{\sigma_k}+k-\sigma_k}~.$$
\end{Proposition}

$\sym$ is a self dual Hopf algebra. It has a pairing (the Hall scalar product)
defined by 
$$\langle h_\lambda, m_\mu \rangle =
\langle s_\lambda, s_\mu \rangle = \delta_{\lambda, \mu}~.$$
An element $f \in \sym$ gives rise to an operator $f^\perp: \sym \to \sym$ 
according to the relation:
\[ \langle fg,h \rangle = \langle g, f^\perp h \rangle \hspace{.1in} \textrm{ for all } g, h \in \sym.\]

We define a ``creation'' operator $\B_m: \sym_n \to \sym_{m+n}$ by: \[\B_m := \sum_{i \geq 0} (-1)^i h_{m+i} e_i^\perp.\]
The following theorem, which states that creation operators construct Schur functions, will become one of the motivations for our new basis of $\Nsym$ (see Definition \ref{def:immaculate}).

\begin{Theorem} (Bernstein \cite[pg 69-70]{Zel}) \label{th:bern}
For all sequences of $\alpha \in \ZZ^m$, 
\[ s_\alpha = \B_{\alpha_1} \B_{\alpha_2} \cdots \B_{\alpha_{m}}(1).\]
\end{Theorem}

\subsection{Non-commutative symmetric functions} 

The algebra $\Nsym$ is a non-commutative analogue of $\sym$ that arises by
considering an algebra with one non-commutative generator at each positive
degree.  In addition to the relationship with the symmetric functions,
this  algebra  has  links  to  Solomon's  descent  algebra  in  type  $A$  \cite{MR},
the  algebra  of  quasi-symmetric  functions  \cite{MR},  and  representation  theory
of  the  type  $A$  Hecke  algebra  at  $q=0$  \cite{KT},  and connections to the theory of combinatorial Hopf algebras \cite{ABS}.  While we will follow the foundational
results  and  definitions  of  references  such  as  \cite{GKLLRT,MR},  we  have  chosen
to use notation here which is suggestive of analogous results in $\sym$.

We define $\Nsym$ as the algebra with generators $\{\HH_1, \HH_2, \dots \}$ and 
no relations. Each generator $H_i$ is defined to be of degree $i$, 
giving $\Nsym$ the structure of a graded algebra. We let $\Nsym_n$ denote the 
graded component of $\Nsym$ of degree $n$. A basis for $\Nsym_n$ are the 
\textit{complete homogeneous functions} 
$\{\HH_\alpha := \HH_{\alpha_1} \HH_{\alpha_2} \cdots \HH_{\alpha_m}\}_{\alpha \vDash n}$ 
indexed by compositions of $n$.  To make this convention consistent, 
some formulas will use expressions that have $H$ indexed by tuples of integers
and we use the convention that $\HH_0=1$ and $\HH_{-r} = 0$ for $r>0$.

There exists a map (sometimes referred to as the forgetful map) which we shall
also denote $\chi: \Nsym \to \sym$ defined by sending the basis element  
$\HH_\alpha$ to the complete homogeneous symmetric function 
\begin{equation}\label{def:chi}
\chi(\HH_\alpha) := h_{\alpha_1} h_{\alpha_2} \cdots h_{\alpha_{\ell(\alpha)}} \in \sym
\end{equation}
and extend this map to all of $\Nsym$ linearly. 

Similar to the study of $\sym$ and the ring of characters for the symmetric groups, the ring of non-commutative symmetric functions of degree $n$ is isomorphic to the Grothendieck ring of  projective representations of the $0$-Hecke algebra. We refer the reader to \cite{KT} for details.
The element of $\Nsym$ which corresponds to the projective representation indexed by $\alpha$ is here denoted $R_\alpha$. The collection of $R_\alpha$ are a basis of $\Nsym$, usually called the \textit{ribbon basis} of $\Nsym$. They are defined through their expansion in the complete homogeneous basis: 
\[ R_\alpha = \sum_{\beta \geq \alpha} (-1)^{\ell(\alpha)-\ell(\beta)} \HH_\beta,
\hskip .1in \textrm{or equivalently,} \hskip .1in \HH_\alpha = \sum_{\beta \geq \alpha} R_\beta. \]
$\Nsym$ has a coproduct structure, which we will not explain in the interest of space.

\subsection{Quasi-symmetric functions}
The algebra of quasi-symmetric functions, $\Qsym$, was introduced in \cite{Ges} 
(see also subsequent
references such as \cite{GR, Sta84}) and this algebra has become a useful tool for algebraic
combinatorics since it is dual to
 $\Nsym$ as a Hopf algebra and contains $\sym$ as a subalgebra. 

As  with  the  algebra  $\Nsym$,  the  graded  component  $\Qsym_n$  is  indexed  by  compositions  of  $n$. 
The algebra is most readily realized within the 
 ring of power series of bounded degree 
$\mathbb{Q}[\![x_1, x_2, \dots]\!]$, and the monomial 
quasi-symmetric function indexed by a composition $\alpha$ is defined as
\begin{equation}
    \label{monomial-qsym}
    M_\alpha = \sum_{i_1 < i_2 < \cdots < i_m} x_{i_1}^{\alpha_1} x_{i_2}^{\alpha_2} \cdots x_{i_m}^{\alpha_m}.
\end{equation}
$\Qsym$ is defined as the algebra with the monomial quasi-symmetric functions as a basis.

We view $\sym$ as a subalgebra of $\Qsym$. In fact, the quasi-symmetric monomial functions
decompose the usual monomial symmetric functions $m_\lambda \in \sym$:
\[ m_\lambda = \sum_{sort(\alpha) = \lambda} M_\alpha.\]

Similar to $\Nsym$, the algebra $\Qsym$ is isomorphic to the Grothendieck ring of finite-dimensional representations of the $0$-Hecke algebra. The irreducible representations of the $0$-Hecke algebra form a basis for this ring, and under this isomorphism the irreducible representation indexed by $\alpha$ is identified with an element of $\Qsym$, the \textit{fundamental quasi-symmetric function}, denoted $F_\alpha$. The $F_\alpha$, for $\alpha \models n$, form a basis of $\Qsym_n$, and are defined by their expansion
in the monomial quasi-symmetric basis: 
\[F_\alpha = \sum_{\beta \leq \alpha} M_\beta.\]
\subsection{Identities relating non-commutative / quasi-symmetric functions} \label{sec:nsymqsymcalc}
The algebras $\Qsym$ and $\Nsym$ form graded dual Hopf algebras. The monomial basis of $\Qsym$ is dual in this context to the complete homogeneous basis of $\Nsym$, and the fundamental basis of $\Qsym$ is dual to the ribbon basis of $\Nsym$.
$\Nsym$ and $\Qsym$ have a pairing $\langle \cdot, \cdot \rangle: \Nsym \times \Qsym \to \mathbb{Q}$, defined under this duality as either $\langle \HH_\alpha, M_\beta \rangle = \delta_{\alpha, \beta}$, or $\langle R_\alpha, F_\beta \rangle = \delta_{\alpha, \beta}$.

We will generalize the operation which is dual to multiplication by a quasi-symmetric function
using this pairing.  For $F, G \in \Qsym$, let $F^\perp$ be the operator which acts on elements
$H \in \Nsym$ according to the relation $\langle H, F G \rangle = \langle F^\perp H, G \rangle$.  

\section{A new basis for $\Nsym$}

We are now ready to introduce our new basis of $\Nsym$.
These functions were discovered while playing with a non-commutative analogue of the 
Jacobi-Trudi identity (see Theorem \ref{thm:JT}). They may also be defined as the 
unique functions in $\Nsym$ which satisfy a right-Pieri rule (see Theorem \ref{thm:Pieri}).

\subsection{Non-commutative immaculate functions}

\begin{Definition}
We define the non-commutative Bernstein operators $\BB_m$ as:
\[ \BB_m = \sum_{i \geq 0} (-1)^i \HH_{m+i} F_{1^i}^\perp~.\]
\end{Definition}

Using the non-commutative Bernstein operators, we can inductively build functions 
using creation operators similar to Bernstein's formula (Theorem \ref{th:bern}) for the Schur functions.

\begin{Remark}
Under the identification of $\sym$ inside $\Qsym$, the generator $e_i$ of $\sym$ is precisely the function $F_{1^i}$ appearing  above.
\end{Remark}

\begin{Definition}\label{def:immaculate}
For any $\alpha = (\alpha_1, \alpha_2, \cdots, \alpha_m) \in \ZZ^m$, 
the \emph{immaculate function} $\fS_\alpha \in \Nsym$ 
is defined as the composition of the operators
$$\fS_\alpha = \BB_{\alpha_1} \BB_{\alpha_2} \cdots \BB_{\alpha_m} (1)~.$$
\end{Definition}

Calculations in the next subsection will show that the elements
$\{ \fS_\alpha \}_{\alpha \models n}$ form a basis for $\Nsym_n$.


\begin{Example}
For $a,b > 0$, $\alpha = (a)$ has only one part, and $\fS_a$ is just the complete homogeneous generator $\HH_a$. If $\alpha = (a, b)$ consists of two parts, then $\fS_{ab} = \BB_a (\HH_b) = \HH_a \HH_b - \HH_{a+1} \HH_{b-1}$.
\end{Example}

\subsection{The right-Pieri rule for the immaculate basis}

\begin{Theorem}\label{thm:Pieri} For a composition $\alpha$, 
the $\fS_\alpha$ satisfy a multiplicity free right-Pieri rule for multiplication by $\HH_s$:
\[\fS_\alpha  \HH_s  =  \sum_{  \alpha  \subset_{s}  \beta}  \fS_\beta.\] 
where the notation $\subset_{s}$ is introduced in Section \ref{sec:compositions}.
\end{Theorem}

\begin{Remark}
Products of the form $\HH_m \fS_\alpha$ do not have as nice an expression as $\fS_\alpha \HH_m$
because they
generally have negative signs in their expansion and there is no obvious containment of
resulting compositions.  For example,
\[ \HH_1 \fS_{13} = \fS_{113} - \fS_{221} - \fS_{32}~. \]
\end{Remark}

\begin{Example}\label{pigeon} The expansion of $\fS_{23}$ multiplied on the right by $\HH_3$ is done below.
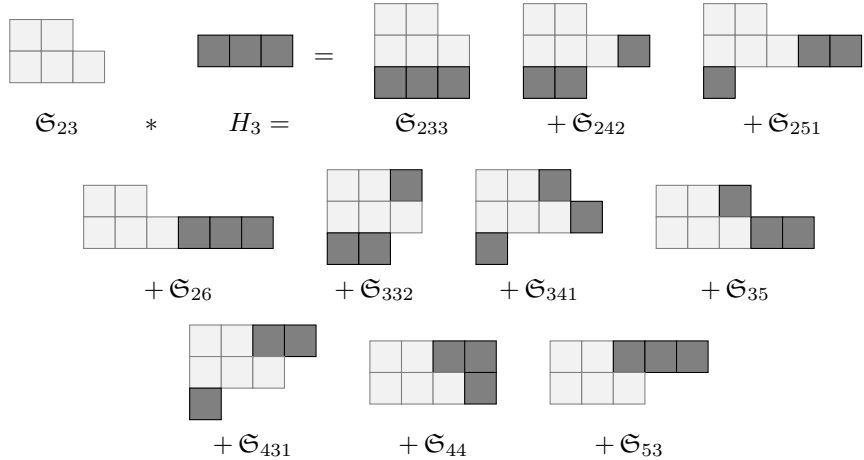
\[
\begin{array}{ccccccccccc}
\tikztableausmall{{X,X},{X,X,X}}& &\tikztableausmall{{\boldentry X,\boldentry X,\boldentry X}}  =   &
\tikztableausmall{{X, X},{X,X, X},{\boldentry X, \boldentry X, \boldentry  X}} &
 \tikztableausmall{{X, X},{X,X, X, \boldentry X},{\boldentry X,\boldentry X}} &
 \tikztableausmall{{X,X},{X,X,X,\boldentry X,\boldentry X},{\boldentry X}}\\
 \fS_{23} &* &\HH_3 = & \fS_{233} &+\, \fS_{242} &+\,\fS_{251}
\end{array}
\]
\[
\begin{array}{ccccc}
& \tikztableausmall{{X,X},{X, X, X, \boldentry X, \boldentry X, \boldentry X}} &
 \tikztableausmall{{X,X,\boldentry X},{X,X,X}, {\boldentry X, \boldentry X}} &
   \tikztableausmall{{X,X,\boldentry X},{X, X, X, \boldentry X}, {\boldentry X}} &
 \tikztableausmall{{X,X,\boldentry X},{X, X, X, \boldentry X, \boldentry X}}\\
&+\,\fS_{26} &+\, \fS_{332}& +\, \fS_{341}& +\, \fS_{35}
\end{array}
\]
\[
\begin{array}{ccccc}
& \tikztableausmall{{X,X,\boldentry X,\boldentry X},{X, X, X}, {\boldentry X}} 
& \tikztableausmall{{X,X,\boldentry X,\boldentry X},{X, X, X, \boldentry X}} 
& \tikztableausmall{{X,X, \boldentry X, \boldentry X, \boldentry X},{X,X,X}}
\\
&+\,\fS_{431} &+ \,\fS_{44} &+\, \fS_{53}
\end{array}
\]
\end{Example}

\subsection{Relationship to the classical bases of $\Nsym$}

We will now develop some relations between the classical bases of $\Nsym$ and the immaculate basis.

\subsubsection{Immaculate tableaux}

\begin{Definition}\label{def:immtab}
Let $\alpha$ and $\beta$ be compositions. An \emph{immaculate tableau} of shape
$\alpha$ and content $\beta$ is a labelling of the boxes of the diagram of
$\alpha$ by positive integers in such a way that:
\begin{enumerate}
\item the number of boxes labelled by $i$ is $\beta_i$;
\item the sequence of entries in each row, from left to right, is weakly increasing;
\item the sequence of entries in the \emph{first} column, from top to bottom,
    is strictly increasing.
\end{enumerate}

An immaculate tableau is said to be \emph{standard} if it has content
$1^{|\alpha|}$.

Let $K_{\alpha, \beta}$ denote the number of immaculate tableaux of shape
$\alpha$ and content $\beta$.
\end{Definition}

We re-iterate that besides the first column, there is no relation on other
columns of an immaculate tableau. Standard immaculate tableaux of size
$n$ are in bijection with set partitions of $\{1, 2, \dots, n\}$ by ordering
the parts in the partition by minimal elements, as was pointed out to us in a discussion with M. Yip.

\begin{Example}\label{ex:immaculatetableau}
There are five immaculate tableau of shape $[4,2,3]$ and content $[3,1,2,3]$: 
\[ \tikztableausmall{{1,1, 1, 3},{2, 3}, {4,4,4}} 
\tikztableausmall{{1,1, 1, 3},{2, 4}, {3,4,4}} 
\tikztableausmall{{1,1, 1, 4},{2,3}, {3,4,4}} 
\tikztableausmall{{1,1, 1, 4},{2, 4}, {3,3,4}} 
\tikztableausmall{{1,1, 1, 2},{3, 3}, {4,4,4}} 
\]
\end{Example}

\subsubsection{Expansion of the homogeneous basis}

\begin{Theorem}\label{thm:Hexpansion}
The complete homogeneous basis $\HH_\alpha$ has a positive, uni-triangular expansion in the immaculate basis. Specifically,
\[ \HH_\beta = \sum_{\alpha \geq_{\mathrm{lex}} \beta} K_{\alpha, \beta} \fS_\alpha,\]
where $K_{\alpha, \beta}$ is the number of immaculate tableaux of shape $\alpha$ and content $\beta$.
\end{Theorem}

\begin{Example}
Continuing from Example \ref{ex:immaculatetableau}, we see that $\HH_{3123} = \cdots + 5 \fS_{423} + \cdots.$
\end{Example}

\begin{Corollary}
The $\{\fS_\alpha : \alpha \vDash n\}$  form a basis of $\Nsym_n$.
\end{Corollary}

\subsubsection{Expansion of the ribbon basis}

We will expand the ribbon functions in the immaculate basis. We first need the notion of a descent.

\begin{Definition} \label{def:descentSIT}
We say that a standard immaculate tableau
$T$ has a descent in position $i$
if  $(i+1)$ is in a row strictly lower than $i$. The descent composition of $T$, $D(T)$, is the composition of the size of $T$ that corresponds to the subset containing all descent positions.
\end{Definition}

\begin{Example}
\label{example:standardization} The standard immaculate tableau below has descents in positions     $\{2, 5, 11\}$.
    The descent composition of $S$ is then $[2,3,6,7]$.
$$ S= \tikztableausmall{{1,2,4,5, 10,11},{3, 6, 7, 8, 9}, {12,13,14,15,16, 17, 18}} $$
\end{Example}

Let $L_{\alpha, \beta}$ denote the number of standard immaculate tableaux of shape $\alpha$ and descent composition $\beta$.

\begin{Theorem}\label{thm:Rpositive}
The ribbon function $R_\beta$ has a positive expansion in the immaculate basis. Specifically 
\[ R_\beta = \sum_{\alpha\geq_\ell\beta} L_{\alpha, \beta} \fS_\alpha.\]
\end{Theorem}

\begin{Example}
There are eight standard immaculate tableaux with descent composition $[2,2,2]$, giving the expansion of $R_{222}$ into the immaculate basis.
$$
\begin{array}{cccccc}
&
\tikztableausmall{{1,2},{3, 4}, {5,6}} &
\tikztableausmall{{1,2},{3, 4,6}, {5}} &
\tikztableausmall{{1,2 , 4},{3}, {5,6}} &
\tikztableausmall{{1,2, 4},{3,6}, {5}}
\tikztableausmall{{1,2, 6},{3,4}, {5}}
\\ R_{222} = &
\fS_{222} &+\, \fS_{231} & +\, \fS_{312} &+\,2\fS_{321}
\end{array}
$$
$$
\begin{array}{cccc}
&
\tikztableausmall{{1,2,4},{3,5,6}}&
\tikztableausmall{{1,2 ,4 ,6},{3}, {5}}&
\tikztableausmall{{1,2, 4, 6},{3,5}}\\
&+\,\fS_{33} &+\, \fS_{411} &+\, \fS_{42}
\end{array}
$$
\end{Example}
\subsection{Jacobi-Trudi rule for $\Nsym$}

Another compelling reason to study the immaculate functions is that they also have an expansion in
the $\HH_\alpha$ basis that makes them a clear analogue of the Jacobi-Trudi rule
of Definition \ref{def:JTformula}.

\begin{Theorem}\label{thm:JT} For a composition
$\alpha = [\alpha_1, \alpha_2, \ldots, \alpha_m]:$
\begin{equation}\label{eq:JTformula}
\fS_\alpha  = \sum_{\sigma \in S_m} (-1)^\sigma \HH_{\alpha_1+\sigma_1 -1, \alpha_2 + \sigma_2 -2, \dots, \alpha_m + \sigma_m - m}.
\end{equation}
\end{Theorem}

\begin{Remark}
This sum is a non-commutative analogue of the determinant of the following matrix:

\[
\begin{bmatrix} 
\HH_{\alpha_1}&\HH_{\alpha_1+1}&\cdots&\HH_{\alpha_1+\ell-1}\\
\HH_{\alpha_2-1}&\HH_{\alpha_2}&\cdots&\HH_{\alpha_2+\ell-2}\\
\vdots&\vdots&\ddots&\vdots\\
\HH_{\alpha_\ell-\ell+1}& \HH_{\alpha_\ell-\ell+2}&\cdots&\HH_{\alpha_\ell}\\
\end{bmatrix}\]
where we have used the convention that $\HH_0 = 1$ and $\HH_{-m} = 0$
for $m>0$.
The non-commutative analogue of the determinant corresponds to expanding this matrix about the first row
and multiplying those elements on the left. 
\end{Remark}

\begin{Remark} One might ask why one would naturally
expand about the first row rather than, say, the first column or the last row.  What
we considered to be the natural analogue of expanding about the first column however
is not a basis; the matrix corresponding to $\alpha = (1,2)$
would be $0$ under this analogue.
\end{Remark}



Of course, the original reason for considering this definition is
the property that they are a lift of the symmetric function corresponding
to the Jacobi-Trudi matrix.

\begin{Corollary}\label{thm:projection}
$\chi( \fS_\alpha) = s_\alpha$.
\end{Corollary}

\subsection{The dual immaculate basis}

Every basis $X_\alpha$ of $\Nsym_n$ gives rise to a basis $Y_\beta$ of $\Qsym_n$ defined by duality; $Y_\beta$ is the unique basis satisfying $  \langle X_\alpha, Y_\beta \rangle = \delta_{\alpha,\beta}$. The dual basis to the immaculate basis of $\Nsym$, denoted $\fS_\alpha^*$, have positive expansions in the monomial and fundamental bases of $\Qsym$. Furthermore, they decompose the usual Schur functions of $\sym$ (see Theorem \ref{thm:decompose}).

\begin{Theorem}\label{thm:MonomialPositive}
The dual immaculate functions $\fS_\alpha^*$ are monomial positive. Specifically they expand as 
\[ \fS_\alpha^* = \sum_{\beta\leq_\ell\alpha} K_{\alpha, \beta} M_\beta.\]
\end{Theorem}

\begin{Theorem}\label{thm:FundamentalPositive}
The dual immaculate functions $\fS_\alpha^*$ are fundamental positive. Specifically they expand as 
\[ \fS_\alpha^* = \sum_{\beta\leq_\ell\alpha} L_{\alpha, \beta} F_\beta.\]
\end{Theorem}

Duality will also yield an explicit expansion of Schur functions into the dual immaculate basis.

\begin{Theorem}\label{thm:decompose}
The Schur function $s_\lambda$, with $\ell(\lambda) = k$ expands into the dual immaculate basis as follows:
\[ s_\lambda = \sum_{\sigma \in S_k} (-1)^\sigma 
\fS_{\lambda_{\sigma_1}+1-\sigma_1, \lambda_{\sigma_2}+2-\sigma_2, \cdots,\lambda_{\sigma_k}+k-\sigma_k}^*\]
where the sum is over permutations $\sigma$ such that $\lambda_{\sigma_i}+i-\sigma_i>0$
for all $i \in \{1,2,\ldots,k\}$.
\end{Theorem}

\begin{Example}
Let $\lambda = (2,2,2,1)$. Then $s_\lambda$ decomposes as:
\[s_{2221} =  \fS^*_{2221} - \fS^*_{1321} - \fS^*_{2131} + \fS^*_{1141},\]
since only the permutations $\sigma \in \{ 1234, 2134, 1324, 2314 \}$ contribute to the sum in the expansion of $s_{2221}$.
There are potentially 24 terms in this sum, but for the partition $(2,2,2,1)$ it is easy to reason that $\sigma_4=4$ and $\sigma_1<3$.
\end{Example}
These combinatorics arise in the paper of Egge, Loehr and Warrington \cite{ELW} when they describe how to obtain a Schur expansion given a quasi-symmetric fundamental expansion. In their language, the terms in this sum correspond to ``special rim hook tableau''.

\subsection{The Littlewood-Richardson rule for immaculate functions}

We prove here that the product $\fS_\alpha \fS_\lambda$ expands positively in the immaculate basis, expanding the notion of a Yamanouchi tableau.
Recall that a \emph{Yamanouchi word} is a word $w$ such that every left prefix of $w$ contains at least as many occurrences of $i$ as $i+1$, for all $i \geq 1$. The content of $w$ is the composition whose $i$-th part is the number of occurrences of $i$.

For partitions $\alpha$ and $\beta$ with $\alpha_i \ge \beta_i$ for all $i$, denote a \emph{skew composition shape} $\alpha // \beta$ by the shape one obtains by superimposing the bottom left boxes of $\alpha$ and $\beta$, and removing the boxes in $\beta$. We denote an \emph{immaculate skew tableau} of shape $\alpha // \beta$ as a filling of this shape, satisfying the rules in Definition \ref{def:immtab}. We denote the \emph{reading word} of a skew immaculate tableau $T$ as the word obtained by reading its entries from right to left in each row, starting from the top row and moving down.

\begin{Theorem}\label{thm:LR}
For a composition $\alpha$ and a partition $\lambda$, the coefficients $c_{\alpha,\lambda}^\beta$ appearing in \[ \fS_\alpha \fS_\lambda = \sum_\beta c_{\alpha, \lambda}^\beta \fS_\beta,\] are non-negative integers. In particular, $c_{\alpha,\lambda}^\beta$ is the number of skew immaculate tableaux of shape $\alpha // \beta$, such that the reading word is a Yamanouchi word of content $\lambda$.
\end{Theorem}


\begin{Example} We give an example with $\alpha = [1,2]$ and $\lambda = [2,1]$.
\[
\begin{array}{ccccccccccccc}
 \tikztableausmall{{\boldentry X},{\boldentry X,\boldentry X}}& &\tikztableausmall{{1,1},{2}}  & =   &
 \tikztableausmall{{\boldentry X},{\boldentry X,\boldentry X},{1,1},{2}} & &
 \tikztableausmall{{\boldentry X},{\boldentry X,\boldentry X,1},{1},{2}} & &
 \tikztableausmall{{\boldentry X},{\boldentry X,\boldentry X,1},{1,2}} & &
 \tikztableausmall{{\boldentry X,1},{\boldentry X,\boldentry X},{1},{2}} & &
 \tikztableausmall{{\boldentry X,1},{\boldentry X,\boldentry X},{1,2}} \\
 \fS_{12} & \hspace{-.2in}*\hspace{-.2in} &\fS_{21} & \hspace{-.2in}=\hspace{-.2in} & \fS_{1221} &\hspace{-.2in}+\hspace{-.2in}& \fS_{1311} &\hspace{-.2in}+\hspace{-.2in}&\,\fS_{132}&\hspace{-.2in}+\hspace{-.2in}&\,\fS_{2211}&\hspace{-.2in}+\hspace{-.2in}&\, \fS_{222}
\end{array}
\]
\[
\begin{array}{ccccccccccccc}
& & & &
 \tikztableausmall{{\boldentry X,1},{\boldentry X,\boldentry X,2},{1}} 
 \tikztableausmall{{\boldentry X,1},{\boldentry X,\boldentry X,1},{2}} & &
 \tikztableausmall{{\boldentry X},{\boldentry X,\boldentry X,1,1},{2}} & &
 \tikztableausmall{{\boldentry X,1},{\boldentry X,\boldentry X,1,2}} & &
 \tikztableausmall{{\boldentry X,1,1},{\boldentry X,\boldentry X,2}} & &
  \tikztableausmall{{\boldentry X,1,1},{\boldentry X,\boldentry X},{2}} \\ 
& & &\hspace{-.2in}+\hspace{-.2in}&\,\ 2\fS_{231} &\hspace{-.2in}+\hspace{-.2in}&\, \fS_{141}&\hspace{-.2in}+\hspace{-.2in}&\, \fS_{24}&\hspace{-.2in}+\hspace{-.2in}&\, \fS_{33}&\hspace{-.2in}+\hspace{-.2in}&\, \fS_{321}
\end{array}
\]
\end{Example}

\subsection{The Murnaghan-Nakayama rule for immaculate functions}

A non-commutative lifting $\Psi_\alpha$ of the power sum basis elements was given in \cite{GKLLRT}. We now state our version of the Murnaghan-Nakayama rule for immaculate functions.

\begin{Theorem}\label{thm:MNrule}
For a composition $\alpha$ and a positive integer $k$,
 \[ \fS_\alpha \Psi_k = \sum_{j=1}^{\ell(\alpha)} \fS_{[\alpha_1, \alpha_2, \dots, \alpha_j+k, \dots, \alpha_{\ell(\alpha)}]} + \sum_{j=0}^{k-1} \fS_{[\alpha, \begin{tiny}\underbrace{0\cdots 0}_j\end{tiny} ,k]}\]
In other words, the sum is over all ways to add $k$ to one of the parts of the composition obtained by padding $\alpha$ with $k$ zeroes at the end.
\end{Theorem}

\begin{Example} One may check that
\[\fS_{132}\Psi_3 = \fS_{432} + \fS_{162}+ \fS_{135}+ \fS_{1323}+ \fS_{13203} + \fS_{132003}.\]
\end{Example}

\begin{Remark}
The nicest form of the Murnaghan-Nakayama rule involves weak compositions (possibly allowing zero as an entry) rather than compositions. There is a signed version of the rule which uses only compositions. In the interest of space, we must omit this rule. It will appear in \cite{BBSSZ3}.
\end{Remark}

\end{document}